\documentclass[12pt]{amsart}

\renewcommand{\epsilon}{\varepsilon}
\renewcommand{\phi}{\varphi}

\newtheorem{Example}{Example}

\begin{document}

\title[Integer Polynomial Maps]{Some Interesting Integer Polynomial Maps}

\author{Alexander Borisov}
\email{borisov@pitt.edu}
\address{Department of Mathematics, University of Pittsburgh, 301 Thackeray Hall, Pittsburgh, PA 15260, USA.}
\thanks{The research of the author was supported in part by the NSA grants H98230-08-1-0129, H98230-06-1-0034 and H98230-11-1-0148.}

\begin{abstract} We introduce three simple polynomial maps with integer coefficients that have interesting dynamical  properties modulo primes.
\end{abstract}

\maketitle

The first map is the simplest, and was discovered while woriking with Mark Sapir on \cite{BS}. 

\begin{Example} (Additive Trap) Define \(F_{at}(x,y)=(u,v)\), where
\[(u,v)=(x^2y,x^2y+xy^2)\]
\end{Example}

Note that modulo any prime \(p\) the \(p\)-th iteration of \(F_{at}\) sends everything to \((0,0)\). Indeed, if \(x\neq 0\), then 
\(\frac{v}{u}=\frac{y}{x}+1\). So after no more than \(p-1\) iterations, we get \(y=0,\) which forces \(x\) and \(y\) to become zero at the next step and forever afterwards. Note that this attractor property of \((0,0)\) does not hold for extensions of \({\mathbb Z}/(p{\mathbb Z})\) or over \(\mathbb Z\) or \(\mathbb C\). Indeed, the map \(F_{at}\) is dominant, so by \cite{BS} its periodic orbits over an algebraic closure of any finite field are Zariski dense.

\begin{Example} (Multiplicative Trap) Define \(F_{mt}(x,y)=(u,v)\), where
\[(u,v)=(x^2y(x-y),2xy^2(x-y))\]
\end{Example}

Note that modulo any prime \(p\) if \(x\neq 0\) and \(\frac{y}{x}\neq 1\), then \(\frac{u}{v}=2\frac{y}{x}\). If \(2\) is a generator of \(({\mathbb Z}/{p\mathbb Z})^*\), then iterations of this map eventually send everything to \((0,0)\).

\begin{Example} (Power Trap) Define \(F_{pt}(x,y)=(u,v)\), where
\[(u,v)=(x^3y(x-y),xy^3(x-y))\]
\end{Example}

Note that modulo any prime \(p\) if \(x\neq 0\) and \(\frac{y}{x}\neq 1\), then \(\frac{u}{v}=(\frac{y}{x})^2\). If \(p=2^k+1,\) then \(k+1\) iterations of this map  send everything to \((0,0)\). For general \(p\) a high enough iteration of \(F_{pt}\) sends to \((0,0)\) all pairs \((x,y)\) with  one of the  coordinates zero or \(\frac{y}{x}\) having 2-primary order in \(({\mathbb Z}/{p\mathbb Z})^*\).

Obviously, with more variables one can create more complicated maps, but it is unclear how much can be actually ``programed" using integer polynomial maps. The following question is very natural in this respect.

{\bf Question.} Does there exist an integer polynomial map with two different fixed points such modulo every prime its sufficiently high iteration will send any initial point to one of the fixed points, depending on whether or not the initial first coordinate is zero?

\end{document}